\documentclass[aps,prl,twocolumn,showpacs,preprintnumbers,amsmath,amssymb,groupedaddress,superscriptaddress]{revtex4-1}

\usepackage{graphicx}
\usepackage{dcolumn}
\usepackage{bm}
\usepackage{subfigure}

\begin{document}

\preprint{APS/123-QED}

\title{Robust convergence in pulse coupled oscillators with delays}

\author{Joel Nishimura}
\affiliation{Center for Applied Mathematics, Cornell University, Ithaca, New York 14853, USA}
\author{Eric J. Friedman}%
\affiliation{Center for Applied Mathematics, Cornell University, Ithaca, New York 14853, USA}
\affiliation{School of ORIE, Cornell University, Ithaca, New York 14853 USA}

\date{\today}

\begin{abstract}
We show that for pulse coupled oscillators a class of phase response curves with both excitation and inhibition exhibit robust convergence to synchrony on arbitrary aperiodic connected graphs with delays.  We describe the  basins of convergence and give explicit bounds on the convergence times.  These results provide new and more robust methods for  synchronization of sensor nets and also have biological implications.

\end{abstract}

\maketitle

Synchronization in systems of pulse coupled oscillators (PCOs) is a fundamental issue in physics, biology and engineering. Examples from nature include synchronization of fireflies \cite{Strogatz}, Josephson junctions \cite{Wiesenfeld}, neurons in the brain \cite{Tateno2007} and the sinoatrial node of the heart \cite{Guevara1982,*Guevara1986}.  In addition to the general study of pulse coupled oscillators \cite{Canavier,*Smeal2010,*Strogatz,*Winfree}, there has been recent interest in their use for synchronization in sensor-networks \cite{Anna,*Konishi,Hu,*Anders,*Fujiwara,*Wang,*KeBo}. However, many of the existing models of PCOs, in particular those for which one can prove analytical results, are limited by strong assumptions.  In this paper we analyze an interesting class of PCOs for which one can prove robust convergence results on arbitrary aperiodic connected graphs, even with propagation delays and a non constant graph topology (such as when spatially embedded nodes are mobile). This class of PCOs is of particular biological relevance because it explicitly includes both inhibition and excitation in the phase response curve (PRC), much like the type II PRCs seen in nature \cite{GoelErmentrout,*Smeal2010,Tateno2007,Guevara1982,*Guevara1986,Anumonwo}.  Additionally, these PCOs provide guidance for the design of engineered systems of PCOs; improving on the current technology, by providing theoretical bounds for robust convergence under propagation delays and covering more diverse topologies. 

Our analysis was motivated by our prior work which used machine learning and genetic algorithms to engineer PRCs which would converge under propagation delays. In that work \cite{inprogress1} we found that such algorithms typically generate a very particular variety of type II PRCs.  As we show below, for engineering applications such as sensor net synchronization \cite{Anna,*Konishi}, these PRCs are superior to those typically used and allow for a precise analysis which appear to differ from the majority of the literature. Namely our analysis does not rely on linear stability, instead our convergence argument makes use of values of the PRC over the entire domain as opposed to derivatives of the PRC at a single point. The analysis also shows that precise normalization of inputs is not required to achieve synchronization with propagation delays, unlike that suggested by the analysis in \cite{Timme,*Flunkert}. 

To begin, we describe the general structure of a PCO model on an arbitrary directed graph under delays. There are $n$ oscillators where oscillator $i$'s state is described by $\phi_i(t) \in [0,1]$.  $\phi_i$ evolves with natural frequency $\dot{\phi} = 1$ and emits a pulse as its phase is reset from $1$ to $0$.  The pulse is received time $\tau <.5$ later by all the successors of $i$, $S(i)$ (predecessors denoted $P(i)$).  Each successor, $j\in S(i)$ adjusts its phase according to its own edge specific PRC, $f_{ji}(\phi_i)$ where $\phi_j\rightarrow \phi_j+f_{ji}(\phi_j)$ (simultaneous signals are processed sequentially in random order).  For example, in the well known model of Mirrolo and Strogatz \cite{Strogatz} this phase adjustment rule is given by $f_{ij}=V^{-1}(\epsilon+V(\phi_i))-\phi_i$ for concave $V$, while the extreme case where  $f_{ij}(\phi_i) = 1-\phi_i$ leads to a resetting of oscillator $i$ with immediate firing while the other extreme case of $f_{ij}(\phi_i) = -\phi_i$ leads to a resetting of oscillator $i$ without firing.  The most well studied models of PCOs are either purely excitatory ($f_{ij}(\phi_i)\geq 0$) or purely inhibitory ($f_{ij}(\phi_i)\leq 0$).

For the sake illustration, consider the PRCs we denote ``strong reseting'' (SR), where for some $B_0\in (\tau,1)$, $f_{ij}(\phi_i)=-\phi_i$ for $0 \le \phi_i \le B_0$ otherwise $f_{ij}=0$. Synchrony is clearly a solution for these curves; every oscillator is simply reset to $0$ time $\tau$ after all oscillators fire. To study this solution consider the time $1+\tau$ map, $H$.

It turns out there is also a clear way to understand part of the basin of the SR system near synchrony. Denote $\phi^+(t)=\mathrm{max}_i (\phi_i(t) )$, $\phi^-(t)=\mathrm{min}_i (\phi_i(t) )$ and $\rho(t)=\phi^+(t)-\phi^-(t)$. Furthermore let $\rho_0(x,y)=\mathrm{min}(x-\tau , 1-y+\tau)$, where $x$ and $y$ are some system parameter. Then in the SR case if at time $t'$ no signals are en route and the range $\rho(t')< \rho_0(B_0,B_0)$ then a careful analysis of the system shows that $$H(\phi_i)=\mathrm{min}( \phi_i+\tau, \mathrm{min}_{j\in P(i)} (\phi_j) ).$$ Notice that if an oscillator $j$ succeeds an oscillator with phases $\phi^-$ then $H(\phi_j)=\phi^-$. In this way the minimum spreads, first to the successors of the minimum and then to the successors of the successors and so on. If the graph is aperiodic, then there exists some $d$ such that $d$ applications of the successor function, denoted: $S^d$, includes the whole graph.  Thus on aperiodic graphs this process leads to synchronization.

SR PRCs harness inhibition to stabilize synchrony, yet in many cases (such as when $\rho>.5$) excitation is useful, as it can allow an oscillator to push forward a cascade. It is possible to augment SR PRCs with excitation while preserving similar convergence bounds and methods of analysis, in the special case where the graph is undirected (directed aperiodic graphs are dealt with by a later theorem). 
Consider the PRC we denote ``strong firing'' (SF) where $f_{ij}(\phi_i)=-\phi_i$ for $0 \le \phi_i \le B_0$ and otherwise $f_{ij}(\phi_i)=1-\phi_i$.  Notice that in this situation oscillators always end up with phase $0$ after they receive a signal, either by being reset or if the phase is greater than $B_0$, by firing. 

To understand the map $H$ in the SF case consider the excitation and inhibition separately, and again constrain $\rho(t')<\rho_0(B_0,B_0)$. It can be verified that between applications of $H$ every oscillator fires. Let $\lambda_i(t_0)$ be the next time that $i$ fires after some time $t_0$.  During an application of $H$ excitation leads to: $$\lambda_i(t_0)=\mathrm{min}(t_0+1-\phi_i(t_0),\mathrm{min}_{j\in P(i)} \lambda_j(t_0)+\tau).$$
Applying this map to an undirected graph implies that all oscillators fire within $\tau$ of their predecessors. The effect of the inhibition then behaves very similar to that of the SR case, where the relative phase differences are captured by $\mathrm{min}(-\lambda_i+\tau, \mathrm{min}_{j\in P(i)} -\lambda_j)$. An interesting feature of the min map in the inhibition is that it preserves the $\tau$ predecessor-successor phase difference, which implies that after a single iteration of $H$ the oscillators will no longer receive signals in the excitatory regime, leading to the same type of convergence as with SR PRCs.  

\begin{figure}
  \includegraphics[width=\linewidth]{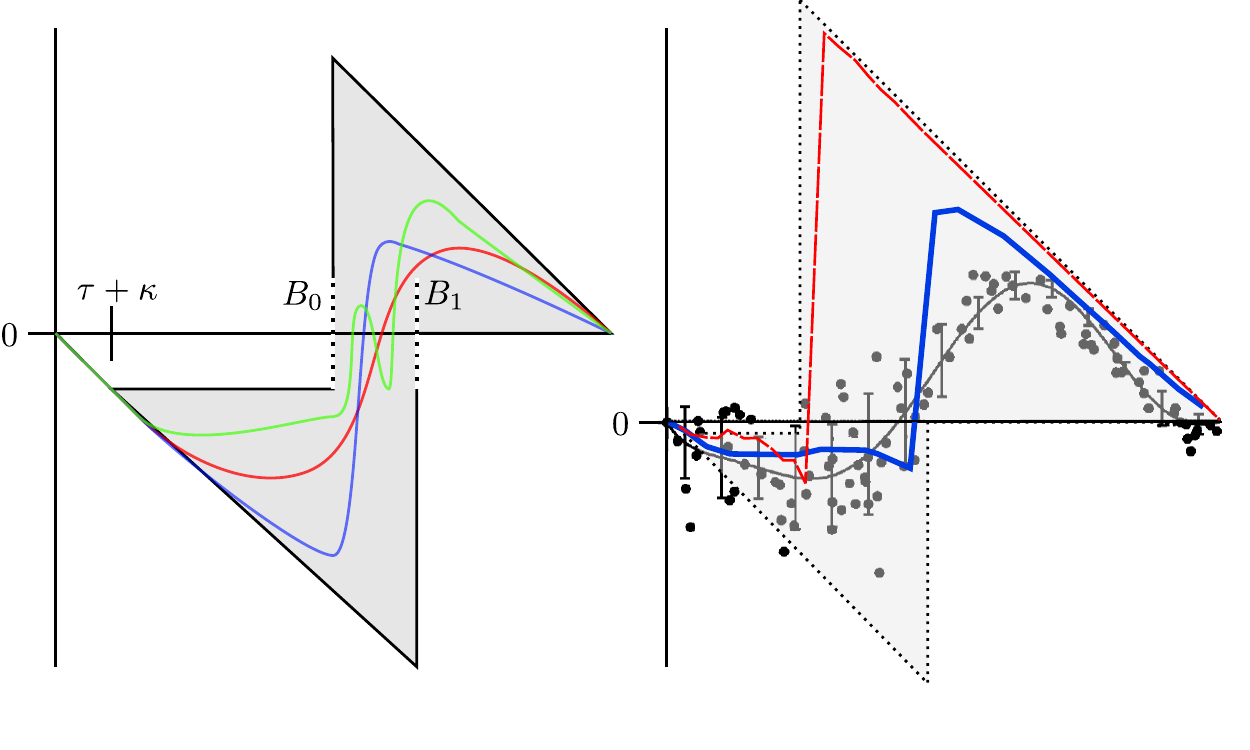}\\
  \caption{(Left) The value of $\tau$ and the shape of S2 PRCs determine the values of $B_0$ and $B_1$ and thus the strength of the convergence results. (Right) Displayed are three different empirically generated PRCs, where the red dashed curve is from ventricular heart cell in an embryonic chick \cite{Guevara1986}, the solid blue curve is from rabbit sinus node cells \cite{Anumonwo} and the data points with the line fit to them is from a low-threshold spiking GABAergic interneurons \cite{Tateno2007}.}\label{Fplot} 
\end{figure}

These two phase response curves converge because of the way a modified min map spreads across a graph, directed aperiodic in the case of SR and undirected aperiodic for SF. While the convergence of SR and SF PRCs is easy to understand other PRCs converge more robustly (as will be demonstrated later) and require a more general argument. Consider the family of PRCs, ``strong type II'' (S2) which have the following requirements. The first requirement has an initial ``reset zone'', $f_{ij}(\phi_i)=-\phi_i$ for $0\leq \phi_i\leq \tau+\kappa$, where $\kappa>0$.  In addition the response curve must be slightly less than $\tau$-inhibitory, in that $f_{ij}(\phi_i)\leq -\tau-\kappa$, on $[\tau,B_0]$ and then be excitatory, $f_{ij}(\phi_i)\ge 0$, for $\phi_i>B_1\in[B_0,1)$. An illustration of S2 curves can be seen in Fig. \ref{Fplot}.  Notice that in $(B_0,B_1)$ there are no restrictions on $f_{ij}$, though the smaller this region is the more general the convergence results will be.  Similarly, it is not necessary that the PRC encode the value of $\tau$ in its shape, only that it is greater than $\tau$ inhibitory.

We can now generalize the results for the SF and SR PRCs to this new class. Furthermore this result remains true even if the graph changes over time. Indeed, suppose that the graph changes every $1+\tau$ time, giving graphs $G_1, G_2,\ldots G_k$. This leads to our main result: If each $G_k$ has no isolated nodes, if there exists $d$ such that $S_l(S_{l+1}(\ldots S_d(v) \dots ))=V$ for an infinite number of $l$, and if the initial range $\rho <\mathrm{min}(B_0-\tau , 1-B_1+\tau)$ the system will converge to synchrony.  Furthermore if the previous conditions hold for all $l$ then convergence occurs before time $t^* = \rho d/ \mathrm{min} (\tau,\kappa)$. The graph conditions might at first seem onerous, but many such examples abound. For example, a grid that suffers a random edge failure each time period would suffice, as would one where each $G_k$ is a different random tree.
While the previous results for SF and SR PRCs were based on an understanding of how one phase spreads across the graph, this argument follows by focusing on the worst case sets of signals an individual oscillator can receive. 

First: if the initial conditions are contained in an interval of size  $\rho$ and there are no signals in transit, it is easy to see that this will remain true under iterations of the time $1+\tau$ map.  To see this, translate time so that the oscillators with the largest phase are just about to fire at time $t=0$: $\phi^+(0)=1^-$.  The key insight is that all signals will occur within at most time $\rho+\tau$, since each oscillator can only fire once in that time, and any oscillator that has not fired will be in the excitatory region. In addition, note that the oscillator with the largest phase can not receive a signal until at least time $\tau$ so will always be inhibited by at least $\tau$. Thus after time $1+\tau$ it will be at most about to fire again, since it was inhibited at least $\tau$ and never excited. A careful analysis of the remaining oscillators using these insights shows that the time $1+\tau$ map does not increase the size of the interval of phases.

Next, we consider oscillators with $\phi_i(0)\leq 1-\epsilon$ where $\epsilon = \min(\kappa,\tau)$. Such an oscillator will not fire until at least time $\epsilon$. Now consider the successors of such an oscillator, $j$ with $\phi_j(0)> 1-\epsilon$.  Oscillator $j$ will receive a signal from $i$ at $t>\epsilon+\tau$ when it is in inhibitory region, so will be inhibited by a sufficient amount such that $\phi_j(1+\tau)\leq 1-\epsilon$.  A careful argument also shows that if $\phi_j(0)\leq 1-\epsilon$ this will remain at time $t+\tau$. Iterating this argument $d$ times (recall that $S^d$ is the complete graph) on the oscillator with the smallest phase, shows that the time $d(1+\tau)$ map will reduce the size of the phase interval by at least $\epsilon$, or if it is less than $\epsilon$ the phases will completely synchronize, proving our convergence result.

\begin{figure}[hbtp]
\includegraphics[width=\linewidth]{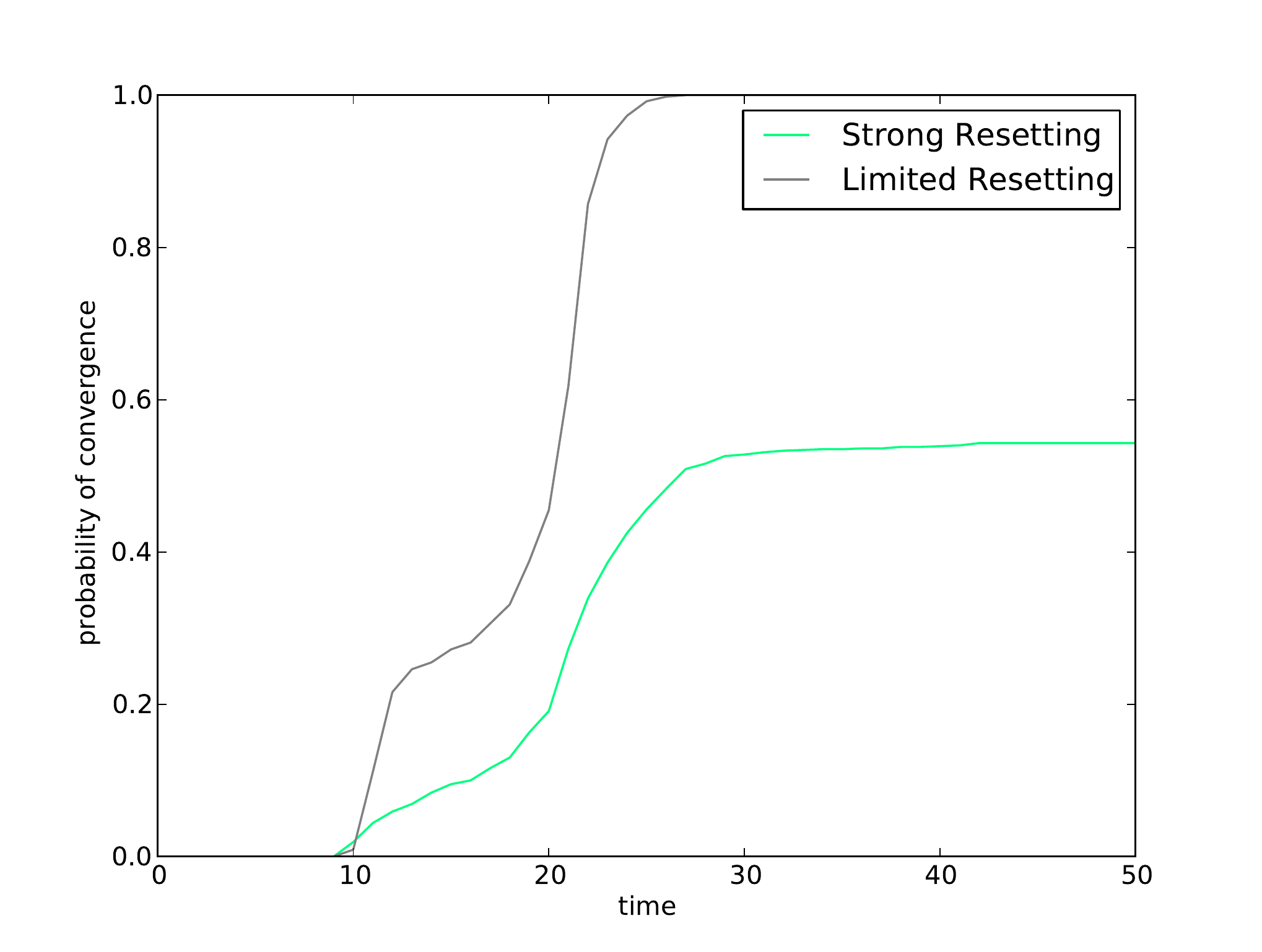}
\caption{Limited Resetting offers more robust convergence on a binary tree with a single triangle and uniform random initial conditions than SF. (The 'S' like shape of the graph is due to whether the system converges to a phase on the same side of the tree as the triangle or far.)  } \label{Basin}
\end{figure}

The benefit of the more general class of phase response curves can be seen in Fig. \ref{Basin}. In Fig. \ref{Basin} the PRC ``Limited Resetting'', which has values strictly greater than $-.1$, is more likely to converge than SR PRCs when run on a slightly modified binary tree with uniform initial conditions. Indeed the combination of inhibition and cuts that divide a graph into multiple disjoint subgraphs can lead to solutions where nodes along the cut never fire (and where the phases are never smaller than the critical $\rho$).  The risk of these non synchronous solutions increases with the amount of inhibition in the PRCs, the number of different disjoint subgraphs and the size of the cuts.

\begin{figure}
\subfigure{\includegraphics[width=2.8in]{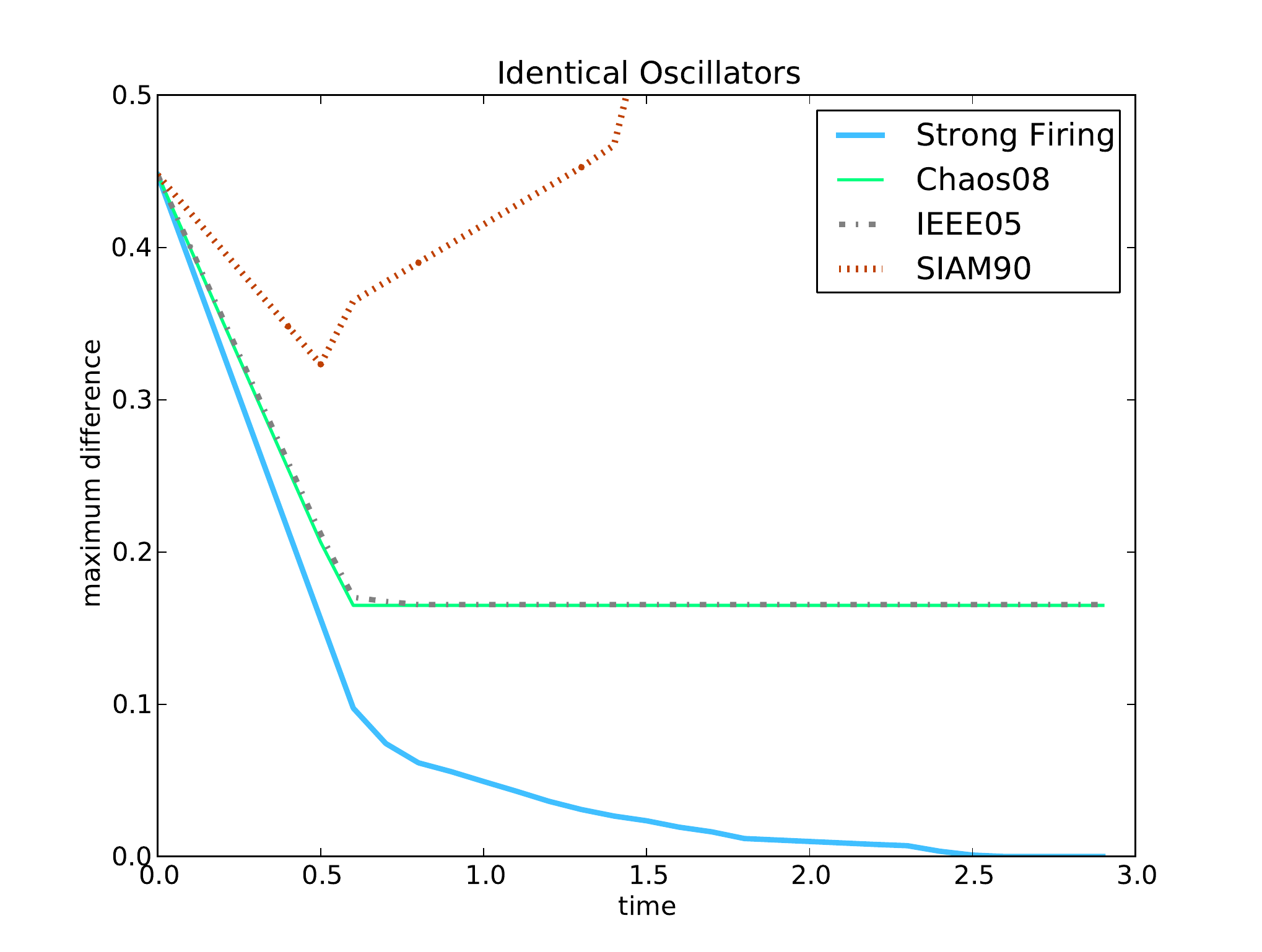}}
\subfigure{\includegraphics[width=2.8in]{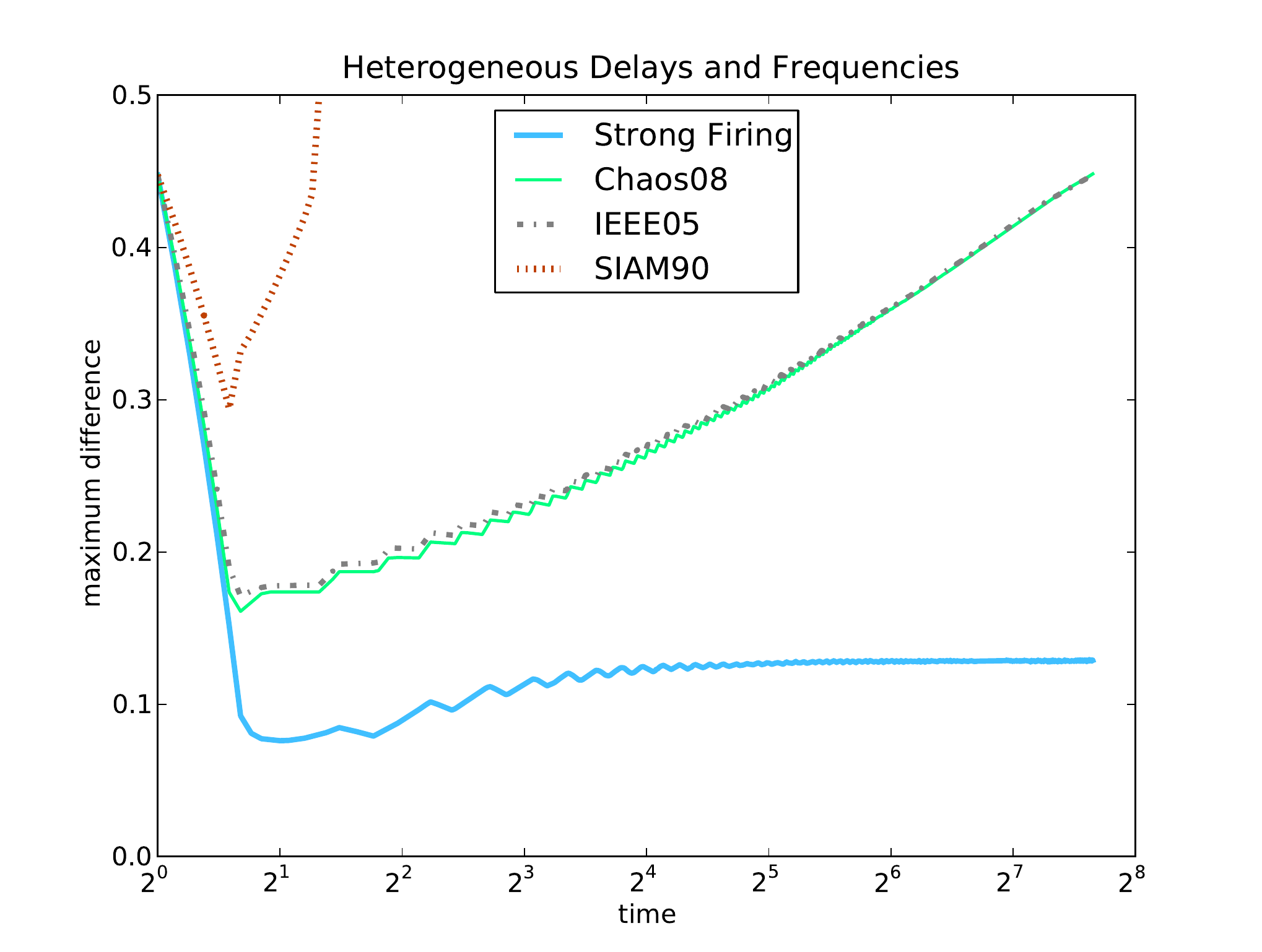}}\\
\subfigure{\includegraphics[width=2.8in]{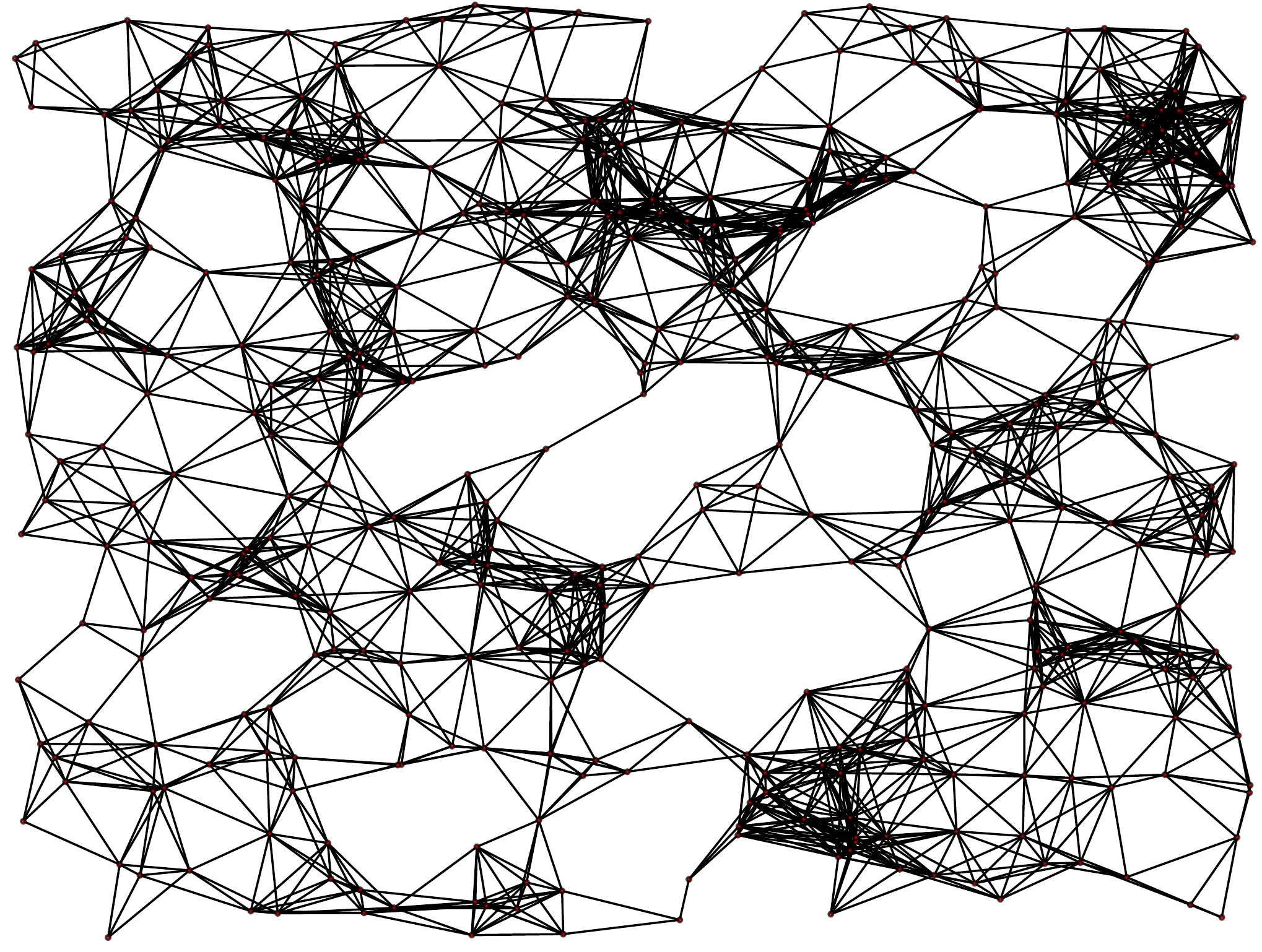}}
\caption{(a) Under the theoretical setting with a uniform delay of $5\%$ S2 curves outperform others from Chaos08 \cite{Konishi}, IEEE05 \cite{Anna}, and SIAM90 \cite{Strogatz} (b) This trend continues under more realistic settings, with frequency error up to $2.5\%$ and transmission and processing delay also up ot $2.5\%$ (c) The random geometric graph on which these simulations was run. \cite{PlotNote} } \label{RGC}
\end{figure}

An important application of this result arises in sensor networks \cite{Anna,*Konishi,Hu,*Anders,*Fujiwara,*Wang,*KeBo}, that is, collections of many small sensors which communicate over radio frequencies. There has been great interest in the use of PCOs to provide a simple and robust mechanism to synchronize sensor networks.  Such systems have intrinsic propagation and processing delays. In addition, most sensor networks have a complex graph structure, and a not necessarily constant network topology. However, most theoretical analyses ignore delays and assume the complete graph. As seen in Fig \ref{RGC}, these assumptions are quite significant. In particular under the assumptions in our first result, the type II PRCs converge rapidly and robustly, while the currently most popular PRCs \cite{Anna,*Konishi} have bounded error but fail to synchronize exactly. In addition, even when we consider more realistic conditions, such as variations in propagation delay times and heterogeneous oscillator frequencies, the type II PRCs still perform well, while the top competitors do not.

Furthermore, the system is relatively robust to error. For example, Notice that the provable basin of attraction is the greatest when $B_0=B_1=.5+\tau$, in which case the the system converges so long as $\rho<.5$.  In this case, with probability $1$ the synchronous solution is robust to any single random error in the oscillator phases. 

Our detailed convergence analysis allows us to ``tune'' the PRCs for specific objectives. For example, as discussed in \cite{Konishi} in sensor-net applications it is  important that the oscillators be allowed to ``sleep'' for as long as possible between firings. In our model this corresponds to choosing a large interval $[B_0,B_1]$ and setting the response curve to $0$ in that region.  However, choosing  a large interval decreases the basin of attraction for the synchronous state. A simple and robust solution for this problem, would be to start with small interval $[B_0,B_1]$ and then expand it over time. For example, one could start with $B_0=B_1=1/2+\tau$  which allows for a large initial interval of phases and then reduce $B_0$ (and increase $B_1$) by $\tau$ or less every $2d(1+\tau)$ units of time, stopping when $B_0=2\tau$. This provides provable convergence with a  long sleep period outside of an initial period.  

In the case where the delay $\tau$ is known in advance then the oscillators can simulate the arrival of their own signal to themselves $\tau$ time after they fire. This modification has the effect of adding self loops to the graphs, vastly increasing the kind of underlying sensor networks that this system performs well on.  Indeed, if ever node has a self loop the only additional requirement on the graph is that the graph never becomes permanently unconnected. 

Type II PRCs have also been seen in many places in nature \cite{GoelErmentrout,Smeal2010}. As seen in Fig. \ref{Fplot}, actual phase response curves taken from cells in the heart and from some cortical interneurons are described by S2 PRC suggesting that the stability of synchronous solutions in these settings may be described by our system. Furthermore, that this family of curves was discovered by a genetic algorithm \cite{inprogress1} lends credence to both the evolability and the performance of such PRC for providing synchrony.  

However, thus far the discussion has included only unweighted graphs, yet many biological graphs are weighted. An extension of this model that allows for weighted graphs is explored in \cite{inprogress2}. This extension allows each edge's PRC to be weighted such that $f_{ij}(\phi)=-\phi$ for $\phi\le w_{ij}$ and less than $w_{ij}$ for $x\in [w_{ij},B_0]$.  In this situation convergence follows from the condition that $\Sigma_j w_{ij}> \tau$.  Thus generalizing the result to weighted graphs.

This weighted version can also be used to give interesting analytic results in the situation where we have some advance knowledge of the underlying graph topology. For example, if the indegree of the graph is known to be at least $k$, then the PRCs only need to be ``resetting'' over the interval from $0$ to $\tau/k$.  This implies that in the limit of high indegree, the resetting region approaches the origin and the only requirement on the magnitude of the inhibition is that the PRC have a slope of $-1$ at the origin and be nonzero in the inhibitory domain. 

One can also use high indegree to provide strong probabilistic convergence results. For example, for any S2 PRC that is nonzero in the excitatory region there exists a constant $c$ such that if the indegree of the connection graph is larger than $c\log(n/\epsilon)$ then under uniform random initial conditions the system will converge to synchrony with probability of at least $\epsilon$.  In general, for a well chosen S2 phase response curve one can give explicit bounds on the probability of convergence for any non zero initial distribution of oscillator phases by analyzing the indegrees.

Other modifications are possible too. One such modification is the addition of ``quiescent'' period $q$, where oscillators ignore signals for time $q$ after receiving a signal. This modification allows for a more general class of PRCs that are allowed more freedom between $0$ and $\tau$.  This new system can be shown to converge on all strongly connected graphs, even periodic graphs.  Another modification is to allow for the oscillator to adjust not only their PRC but also their frequency.  Numerical results suggest that the proper choice of the frequency response curve allows for oscillators with a more robust region of convergence. 

In summary, the family of S2 phase response curves that was introduced is relatively general and includes curves that were empirically found in systems which synchronize in nature. Furthermore, we outlined a proof that this family of PCOs has a robust region in which it converges to synchrony on strongly connected aperiodic graphs. This convergence remains even in the presence of uniform time delay and particular mutations in the graph. It was then noted that in numerical trials this method outperforms similar PCO based time synchronization methods for wireless sensor networks. This advantage remained in numerical runs even when heterogeneous time delays, frequencies, and random errors were introduced.

We thank Bard Ermentrout, Jon Kleinberg, David Shmoys, Steven Strogatz, and  Xiao Wang
for useful conversations and the EECS at UC Berkeley for its hospitality. This research has been supported in part by the NSF under grant CDI-0835706.

\bibliography{OscReferences}

\begin{thebibliography}{22}
\bibitem [{Plo()}]{PlotNote}%
  \BibitemOpen
  \href@noop {} {}\bibinfo {note} {Data points are taken when oscillators are
  within a range of .5 and no signals are en route}\BibitemShut {NoStop}%
\end{thebibliography}

\end{document}